\documentclass[12pt]{amsart}
\usepackage{graphicx,enumerate,amsmath,amssymb,afterpage,cite}
\usepackage[sc]{mathpazo}

\newcommand{\dif}{\mathrm{d}}

\newcommand{\parderiv}[2]{\frac{\partial #1}{\partial #2}}
\newcommand{\deriv}[2]{\frac{\mathrm{d} #1}{\mathrm{d} #2}}
\newcommand\BS\boldsymbol
\newcommand\sym{\mathrm{\,Sym}}
\newcommand\picturehere[1]{\includegraphics[width=0.5\textwidth]{#1}}

\numberwithin{equation}{section}
\numberwithin{figure}{section}
\numberwithin{table}{section}

\begin{document}

\title[A simple stochastic parameterization for reduced models]
{A simple stochastic parameterization for reduced models of multiscale dynamics}

\author{Rafail V. Abramov}

\address{Department of Mathematics, Statistics and Computer
  Science\\University of Illinois at Chicago\\851 S. Morgan st. (M/C
  249)\\ Chicago, IL 60607}

\email{abramov@math.uic.edu}


\date{\today}

\pagestyle{myheadings}

\begin{abstract}
Multiscale dynamics are frequently present in real-world processes,
such as the atmosphere-ocean and climate science. Because of time
scale separation between a small set of slowly evolving variables and
much larger set of rapidly changing variables, direct numerical
simulations of such systems are difficult to carry out due to many
dynamical variables and the need for an extremely small time
discretization step to resolve fast dynamics. One of the common
remedies for that is to approximate a multiscale dynamical systems by
a closed approximate model for slow variables alone, which reduces the
total effective dimension of the phase space of dynamics, as well as
allows for a longer time discretization step. Recently we developed a
new method for constructing a deterministic reduced model of
multiscale dynamics where coupling terms were parameterized via the
Fluctuation-Dissipation theorem. In this work we further improve this
previously developed method for deterministic reduced models of
multiscale dynamics by introducing a new method for parameterizing
slow-fast interactions through additive stochastic noise in a
systematic fashion. For the two-scale Lorenz 96 system with linear
coupling, we demonstrate that the new method is able to recover
additional features of multiscale dynamics in a stochastically forced
reduced model, which the previously developed deterministic method
could not reproduce.
\end{abstract}

\maketitle

\section{Introduction}

Multiscale dynamics are common in applications of contemporary
science, such as geophysical science and climate change prediction
\cite{FraMajVan,Has,BuiMilPal,Pal3,BraBer,MajFraCro,KraKonGhi}.
Multiscale dynamics are typically characterized by the time and space
scale separation of patterns of motion, with fewer slowly evolving
variables and much larger set of faster evolving variables. This
time-space scale separation often made direct numerical computation of
the dynamics on the slow time scale quite difficult in real-world
applications, which led to the development of multiscale computational
methods \cite{ELiuVan,FatVan}. These methods make use of the averaging
formalism \cite{Pap,Van,Vol} to allow for large time discretization
steps for the computation of the slow part of the dynamics. However,
in very large systems with many fast variables even these methods are
computationally expensive.

As a different alternative to direct numerical simulation of the
complete multiscale model with all variables, it has long been
recognized that, if a closed simplified model for the slow variables
alone was available, one could use this closed slow-variable model
instead to simulate the statistics of the slow variables. In order to
derive such a reduced model, one usually represents the process as a
general two-scale dynamical system of the form
\begin{equation}
\label{eq:dyn_sys}
\deriv{\BS x}t=\BS F(\BS x,\BS y),\qquad
\deriv{\BS y}t=\frac 1\varepsilon\BS G(\BS x,\BS y),
\end{equation}
where $\BS x\in\mathbb R^{N_x}$ and $\BS y\in\mathbb R^{N_y}$ are the
state vectors, and $\BS F$ and $\BS G$ are nonlinear differentiable
functions. The scaling parameter $\varepsilon\ll 1$ is used to
separate the time scales in \eqref{eq:dyn_sys} into slow (that is,
$\BS x$), and fast (that is, $\BS y$). The scale separation parameter
$\varepsilon$ above is introduced for convenience of presentation,
because, as shown below, the method developed here does not require
its explicit presence to function. Besides, often real-world processes
do not have any distinct time-scale separation parameters, and the
fast and slow variables are known empirically from observations.

Under the assumption of ``infinitely fast'' $\BS y$-variables, one
writes the reduced averaged system for slow variables alone as
\begin{equation}
\label{eq:dyn_sys_slow_limiting_x}
\deriv{\BS{\bar x}}t=\BS{\bar F}(\BS{\bar x}),\qquad\BS{\bar
  F}(\BS{\bar x})=\int_{\mathbb R^{N_y}}\BS F(\BS{\bar x},\BS y)
\dif\mu_{\BS{\bar x}}(\BS y),
\end{equation}
where $\mu_{\BS{\bar x}}$ is the invariant distribution measure of the
uncoupled fast dynamics where $\BS{\bar x}$ is a fixed parameter:
\begin{equation}
\label{eq:dyn_sys_fast_limiting_z}
\deriv{\BS z}t=\BS G(\BS{\bar x},\BS z).
\end{equation}
Under the assumption of ergodic $\mu_{\BS{\bar x}}$, the measure
integral in \eqref{eq:dyn_sys_slow_limiting_x} can be replaced with
the time average along a single trajectory $\BS z_{\BS{\bar x}}(t)$ of
\eqref{eq:dyn_sys_fast_limiting_z} for specified parameter $\BS{\bar x}$:
\begin{equation}
\label{eq:F_average}
\BS{\bar F}(\BS{\bar x})=\lim_{T\to\infty}\frac 1T\int_0^T\BS F(\BS{\bar x},
\BS z_{\BS{\bar x}}(t))\dif t.
\end{equation}
From \eqref{eq:dyn_sys_slow_limiting_x}--\eqref{eq:F_average}, it is
clear that the computation of the averaged function $\BS{\bar
  F}(\BS{\bar x})$ is not a simple task; in fact, it is rarely
available explicitly (except for some special cases where either the
invariant distribution measure $\mu_{\BS{\bar x}}$ or the solution
$\BS z_{\BS{\bar x}}(t)$ are known as explicit formulas). A good
example of the system where the invariant distribution measure is
known explicitly is the Ornstein-Uhlenbeck process \cite{OrnUhl}; the
invariant distribution measure of the Ornstein-Uhlenbeck process is a
Gaussian distribution with explicitly known mean state and covariance
matrix. Because of this convenient property, the Ornstein-Uhlenbeck
process is popular in the area of stochastic modeling of geophysical
and other real-world processes.

The ``upgrade'' from the deterministic reduced model in
\eqref{eq:dyn_sys_slow_limiting_x} is a stochastic model of the form
\begin{equation}
\label{eq:dyn_sys_slow_limiting_stoch}
\dif\BS{\bar x}=\BS{\bar F}(\BS{\bar x})\dif t+\BS\sigma(\BS{\bar
  x})\dif\BS W_t,
\end{equation}
where $\BS W_t$ is a Wiener process, and $\BS\sigma(\BS{\bar x})$ is a
matrix. The stochastic reduced model of the form
\eqref{eq:dyn_sys_slow_limiting_stoch} is generally considered to be
more advantageous to the deterministic reduced model in
\eqref{eq:dyn_sys_slow_limiting_x}, because the random noise can be
used to parameterize ``noise-like'' influence from fast variables onto
slow dynamics, which is entirely lacking in the deterministic reduced
model in \eqref{eq:dyn_sys_slow_limiting_x}. Sometimes, deterministic
reduced models of atmospheric processes fail to adequately represent
atmospheric variability; it is thought that a significant portion of
atmospheric variability is essentially noise-like, and fitting it with
deterministic forcing terms does not seem to be a good way to produce
an adequate approximation. Additionally, the deterministic slow
dynamics in \eqref{eq:dyn_sys_slow_limiting_x} tend to produce the
attractor of the slow dynamics with significantly lower dimension than
that for the full two-scale dynamics in \eqref{eq:dyn_sys}, especially
for weakly chaotic slow variables. At the same time, if the stochastic
forcing in \eqref{eq:dyn_sys_slow_limiting_stoch} is strong enough, it
should ``diffuse'' the invariant manifolds of
\eqref{eq:dyn_sys_slow_limiting_x}, thus inflating the dimension of
the set containing limiting dynamics.

Nonetheless, in many applications, the averaged dynamics of the
deterministic type \eqref{eq:dyn_sys_slow_limiting_x} or stochastic
type \eqref{eq:dyn_sys_slow_limiting_stoch} are not available
explicitly. As a result, numerous approximate closure schemes were
developed for multiscale dynamical systems
\cite{CroVan,FatVan,MajTimVan,MajTimVan2,MajTimVan3,MajTimVan4}, which
are all based on the averaging principle over the fast variables
\cite{Pap,Van,Vol}. Some of the methods (such as those in
\cite{MajTimVan,MajTimVan2,MajTimVan3,MajTimVan4}) replace the fast
nonlinear dynamics with suitable stochastic processes \cite{Wilks},
discontinuous Markov jump processes \cite{KatVla}, or conditional
Markov chains \cite{CroVan}, while others \cite{FatVan} provide direct
closure by suitable tabulation and curve fitting. Reduced stochastic
dynamics were used to model global circulation patterns
\cite{FraMaj,Bra,NewSarPen,WhiSar,ZhaHel}, and large-scale features of
tropical convection \cite{MajKho,KhoMajKat}. However, it seems that
all these approaches require either extensive computations to produce
a closed model (for example,\cite{CroVan,FatVan} require multiple
simulations of fast variables alone with different fixed states of
slow variables), or somewhat {\em ad hoc} determination of closure
coefficients by matching areas under the time correlation functions
\cite{MajTimVan,MajTimVan2,MajTimVan3,MajTimVan4}. Another interesting
method was recently developed in \cite{AzeBerTim}, however, again,
somewhat {\em ad hoc} approach was used to compute the closure
(namely, slow variables of a multiscale system were treated as if
their dynamics were generated by a running average of an
Ornstein-Uhlenbeck process).

In the two recent works \cite{Abr9,Abr10} the author developed a
relatively simple and straightforward method of constructing the
deterministic reduced model for slow variables of a multiscale model
with nonlinear and multiplicative coupling, which required only a
single computation of certain statistics of the fast dynamics with a
fixed state of the slow variables, located in the region where the
slow dynamics usually evolve. The method was based on the first-order
Taylor expansion of the averaged coupling term with respect to the
slow variables, which was computed using the Fluctuation-Dissipation
theorem
\cite{Abr5,Abr6,Abr7,AbrMaj4,AbrMaj5,AbrMaj6,MajAbrGro,Ris}. It was
demonstrated through the computations with the appropriately rescaled
two-scale Lorenz 96 model \cite{Lor,LorEma} that, with nonlinear and
multiplicative coupling in both slow and fast variables, the developed
reduced model produced good approximation to the statistics of the
full two-scale Lorenz 96 model. Among the advantages of the developed
method were its simplicity and explicit formulation. Additionally,
existing zero-order models of this kind for the Earth's atmosphere
(such as the T21 barotropic model \cite{AbrMaj6,Fra,Sel}) can be
retrofitted with the new deterministic correction term emerging from
the theory in \cite{Abr9}.

In the current work, we introduce a method for computing a consistent
approximation of the form \eqref{eq:dyn_sys_slow_limiting_stoch} to
the multiscale dynamics in \eqref{eq:dyn_sys}. The method is aimed at
stochastic parameterization of general complex nonlinear multiscale
dynamics with many variables, and has essentially the same
implementation restrictions as the method for deterministic reduced
models we developed previously in \cite{Abr9,Abr10}. We test the new
method on the two-scale Lorenz 96 model, where only the linear part of
the coupling is enabled for the simplicity of presentation. We
demonstrate through direct numerical simulations that the new method
generally improves properties of statistics of the reduced model, and,
in particular, the injection of random noise into a deterministic
reduced model can make it both more or less chaotic and mixing,
depending on the difference between the dynamical regimes of the
reduced and full two-scale models.

The manuscript is organized as follows. In Section \ref{sec:theory} we
present the general description of the new method, following the
homogenization theory of \cite{PavStu}. In Section \ref{sec:practical}
we lay out the step-by-step computational implementation of the new
method for a general two-scale dynamical system with linear coupling
between the slow and fast variables, which does not contain any
explicit time-scale separation parameters. In Section
\ref{sec:numerical} we test the new method on the two-scale Lorenz 96
model in a range of dynamical regimes with varying chaos, mixing, and
time scale separation between the slow and fast variables. Section
\ref{sec:conclusions} summarizes the results of this work.

\section{General description of the method}
\label{sec:theory}

Here we introduce a new method for the stochastic correction of the
form \eqref{eq:dyn_sys_slow_limiting_stoch} for the deterministic
reduced model \eqref{eq:dyn_sys_slow_limiting_x}. To derive the
method, we use the theoretical framework similar to that applied to
homogenization problems in \cite{PavStu}. To simplify presentation, we
assume that the deterministic averaged function $\BS{\bar F}(\BS{\bar
  x})$ from \eqref{eq:dyn_sys_slow_limiting_x} is already available,
either as an explicit formula, or as an approximation we developed
previously in \cite{Abr9,Abr10}. Below, $\BS x$ is used to denote the
state of the slow variables from the multiscale system
\eqref{eq:dyn_sys}, while $\BS{\bar x}$ denotes the state of the slow
variables from the deterministic reduced model
\eqref{eq:dyn_sys_slow_limiting_x}. The difference between $\BS x$ and
$\BS{\bar x}$ is denoted as $\BS q=\BS x-\BS{\bar x}$. Then, one can
rewrite the multiscale system in \eqref{eq:dyn_sys} in the new
variables as
\begin{subequations}
\label{eq:dyn_sys_q}
\begin{equation}
\deriv{\BS{\bar x}}t=\BS{\bar F}(\BS{\bar x}),
\end{equation}
\begin{equation}
\label{eq:q}
\deriv{\BS q}t=\BS F(\BS{\bar x}+\BS q,\BS y)-\BS{\bar F}(\BS{\bar x}),
\end{equation}
\begin{equation}
\deriv{\BS y}t=\frac 1\varepsilon\BS G(\BS{\bar x}+\BS q,\BS y).
\end{equation}
\end{subequations}
What we see above is that the first equation is already a closed
system from \eqref{eq:dyn_sys_slow_limiting_x}, and given a suitable
approximation for $\BS{\bar F}$ from \cite{Abr9,Abr10}, it can be
solved on its own. Thus, $\BS{\bar x}(t)$ can be treated as a given
function of time. This, in effect, leaves $\BS q$ and $\BS y$ as the
unknown variables, and the first equation in \eqref{eq:dyn_sys_q} can
be dropped. Now, the idea is to apply the averaging formalism to $\BS
q$, obtaining $\BS{\bar q}$ as the next order correction to $\BS{\bar
  x}$. However, the straightforward application of what was done in
\eqref{eq:dyn_sys_slow_limiting_x} to \eqref{eq:dyn_sys_q} leads to
$\BS{\bar q}$ being identically zero for all times as long as its
starting value is zero. In this situation, the reduced model for
$\BS{\bar q}$ should be derived as the It\^o diffusion process of the
corresponding backward Kolmogorov equation restricted to slow time
scale (see \cite{PavStu} for details), under the condition that $\BS
q$ (and, therefore, $\dif\BS q/\dif t$) in \eqref{eq:dyn_sys_q} is
$O(\varepsilon)$. To do that, first we bring the time scale of $\BS q$
in \eqref{eq:dyn_sys_q} to $O(1)$ by rescaling the time $t$ as
$\tau=\varepsilon t$. For the rescaled time $\tau$, from
\eqref{eq:dyn_sys_q} we obtain
\begin{subequations}
\label{eq:dyn_sys_q_rescaled}
\begin{equation}
\deriv{\BS q}\tau=\frac 1\varepsilon\left[\BS F(\BS{\bar x}+\BS q,\BS
  y)-\BS{\bar F}(\BS{\bar x})\right],
\end{equation}
\begin{equation}
\deriv{\BS y}\tau=\frac 1{\varepsilon^2}\BS G(\BS{\bar x}+\BS q,\BS y).
\end{equation}
\end{subequations}
Above, the evolution equation for $\BS{\bar x}$ is no longer needed,
as $\BS{\bar x}(\tau)$ is a given function of rescaled time
$\tau$. Now we write the backward Kolmogorov equation for
\eqref{eq:dyn_sys_q_rescaled}. For that, let
$v(\tau,\tau^\prime,\BS{\bar x},\BS q,\BS y)$, $\tau^\prime\geq\tau$,
be the value of a test function $h(\BS q(\tau^\prime),\BS
y(\tau^\prime))$, given $\BS q(\tau)=\BS q$, and $\BS y(\tau)=\BS
y$. Then, $v(\tau,\tau^\prime,\BS{\bar x}, \BS q,\BS y)$ obeys the
following backward Kolmogorov equation (for reference, see, for
example, \cite{GikSko}):
\begin{equation}
\label{eq:kolmo}
\parderiv{v(\tau,\tau^\prime,\BS{\bar x},\BS q,\BS y)}\tau=\left[\frac
  1{\varepsilon^2} \BS G(\BS{\bar x}+\BS q,\BS y)\cdot\nabla_{\BS
    y}+\frac 1\varepsilon \big[\BS F(\BS{\bar x}+\BS q,\BS y)-\BS{\bar
      F}(\BS{\bar x})\big] \cdot\nabla_{\BS
    q}\right]v(\tau,\tau^\prime,\BS{\bar x},\BS q,\BS y).
\end{equation}
Below, we drop the $\tau^\prime$-dependence from $v$ as it is of no
consequence to what is presented. Observe that the terms in the
backward Kolmogorov equation above are multiplied by different powers
of $\varepsilon$, which leads to the perturbation expansion of the
solution $v(\tau,\BS{\bar x},\BS q,\BS y)$ in powers of $\varepsilon$, and
derivation of the backward Kolmogorov equation for the term which has
the lowest power of $\varepsilon$ in the closed form. Then, its
corresponding It\^o diffusion process will be the evolution equation
for $\dif\BS{\bar q}/\dif t$, as long as $\BS{\bar q}$ is small
enough.

The expansion of $v(\tau,\BS{\bar x},\BS q,\BS y)$ in powers of
$\varepsilon$ is
\begin{equation}
v=v_0+\varepsilon v_1+\varepsilon^2 v_2+\ldots.
\end{equation}
Plugging the expansion above back into the Kolmogorov equation in
\eqref{eq:kolmo} and collecting the terms with matching powers of
$\varepsilon$, we obtain the following relations for each power of
$\varepsilon$:
\begin{subequations}
\begin{equation}
\label{eq:v0}
\BS G(\BS{\bar x}+\BS q,\BS y)\cdot\nabla_{\BS y}v_0=0\quad\mbox{for}\quad\varepsilon^{-2},
\end{equation}
\begin{equation}
\label{eq:v1}
\BS G(\BS{\bar x}+\BS q,\BS y)\cdot\nabla_{\BS y}v_1=-\big[\BS F(\BS{\bar x}+\BS q,\BS y)
-\BS{\bar F}(\BS{\bar x})\big]\cdot\nabla_{\BS q}v_0\quad\mbox{for}\quad\varepsilon^{-1},
\end{equation}
\begin{equation}
\label{eq:v2}
\parderiv{v_0}\tau=\BS G(\BS{\bar x}+\BS q,\BS y)\cdot\nabla_{\BS y}v_2+
\big[\BS F(\BS{\bar x}+\BS q,\BS y)-\BS{\bar F}(\BS{\bar x})\big]\cdot
\nabla_{\BS q}v_1\quad\mbox{for}\quad\varepsilon^0.
\end{equation}
\end{subequations}
Below we consider each of the relations above separately.
\begin{itemize}
\item {\bf Order $\varepsilon^{-2}$.} From the relation in
  \eqref{eq:v0} we determine that $v_0(\tau,\BS{\bar x},\BS q,\BS
  y)=v_0(\tau,\BS{\bar x},\BS q)$, that is, $v_0$ does not depend on
  $\BS y$.
\item {\bf Order $\varepsilon^{-1}$}. Here we use the relation in
  \eqref{eq:v1} to express $v_1$ in terms of $v_0$. We denote the
  flow, generated by \eqref{eq:dyn_sys_fast_limiting_z}, by
  $\phi_{\BS{\bar x}}^s$, so that $\phi_{\BS{\bar x}}^s\BS y$ is the
  solution of \eqref{eq:dyn_sys_fast_limiting_z} forward in time $s$
  with the initial condition $\BS y$, with the obvious identity
\begin{equation}
\label{eq:flow_identity}
\parderiv{}s\phi_{\BS{\bar x}+\BS q}^s\BS y=\BS G(\BS{\bar x}+\BS q,\phi_{\BS{\bar x}+\BS q}^s\BS y).
\end{equation}
Now, consider the integral
\begin{equation}
\label{eq:two_identities}
\begin{split}
u(s,\BS{\bar x},\BS q,\BS y)=\int_s^\infty\big[\BS F(\BS{\bar x}+\BS q,
\phi_{\BS{\bar x}+\BS q}^r\BS y)-\BS{\bar F}(\BS{\bar x})\big]\dif r\cdot
\nabla_{\BS q}v_0=\\=\int_0^\infty\big[\BS F(\BS{\bar x}+\BS q,
\phi_{\BS{\bar x}+\BS q}^r(\phi_{\BS{\bar x}+\BS q}^s\BS y))-\BS{\bar F}
(\BS{\bar x})\big]\dif r\cdot\nabla_{\BS q}v_0,
\end{split}
\end{equation}
where the group property of $\phi_{\BS{\bar x}}^s$ is used in the second equality.
Then, from the first identity in \eqref{eq:two_identities}, it follows
that
\begin{equation}
\label{eq:initial_condition}
\left.\parderiv{}s u(s,\BS{\bar x},\BS q,\BS y)\right|_{s=0}=
-\big[\BS F(\BS{\bar x}+\BS q,\BS y)-\BS{\bar F}(\BS{\bar x})\big]
\cdot\nabla_{\BS q}v_0,
\end{equation}
and from the second identity in \eqref{eq:two_identities} it follows
that $u(s,\BS{\bar x},\BS q,\BS y)$ is in fact an explicit function of
$\phi_{\BS{\bar x}+\BS q}^s\BS y$, that is, $u(s,\BS{\bar x},\BS q,\BS
y)\equiv U(\BS{\bar x},\BS q,\phi_{\BS{\bar x}+\BS q}^s\BS
y)$. However, any $U(\BS{\bar x},\BS q,\phi_{\BS{\bar x}+\BS q}^s\BS
y)$ must obey the transport equation
\begin{equation}
\label{eq:transport}
\begin{split}
\parderiv{}s U(\BS{\bar x},\BS q,\phi_{\BS{\bar x}+\BS q}^s\BS y) =
\nabla U(\BS{\bar x},\BS q, \phi_{\BS{\bar x}+\BS q}^s\BS y)\cdot\parderiv{}s
\phi_{\BS{\bar x}+\BS q}^s\BS y=\\=\BS G(\BS{\bar x}+\BS q,
\phi_{\BS{\bar x}+\BS q}^s\BS y)\cdot\nabla U(\BS{\bar x},\BS q,
\phi_{\BS{\bar x}+\BS q}^s\BS y),
\end{split}
\end{equation}
where the second identity is due to \eqref{eq:flow_identity}, and
which holds for any $s$ including $s=0$. Combining
\eqref{eq:initial_condition} and \eqref{eq:transport} at $s=0$, we
obtain
\begin{equation}
\BS G(\BS{\bar x}+\BS q,\BS y)\cdot\nabla_{\BS y}u(0,\BS{\bar x},\BS q,
\BS y)=-\big[\BS F(\BS{\bar x}+\BS q,\BS y)-\BS{\bar F}(\BS{\bar x})\big]
\cdot\nabla_{\BS q}v_0.
\end{equation}
From the comparison with \eqref{eq:v1} it follows that
$v_1=u(0,\BS{\bar x},\BS q,\BS y)$, that is,
\begin{equation}
\label{eq:v11}
v_1=\int_0^\infty\big[\BS F(\BS{\bar x}+\BS q,\phi_{\BS{\bar x}+\BS q}^s
\BS y)-\BS{\bar F}(\BS{\bar x})\big]\dif s\cdot\nabla_{\BS q}v_0.
\end{equation}
\item {\bf Order $\varepsilon^0$.} Here observe that $v_0$ does not
  depend on $\BS y$, as pointed out above. This means that the average
  of $v_0$ with respect to the invariant distribution measure
  $\mu_{\BS{\bar x}+\BS q}$ of \eqref{eq:dyn_sys_fast_limiting_z} is
  the identity operation, and the same holds for its
  $\tau$-derivative. Then, averaging out \eqref{eq:v2} with respect to
  $\mu_{\BS{\bar x}+\BS q}$ yields,
\begin{equation}
\label{eq:v22}
\begin{split}
\parderiv{v_0}\tau&=\int_{\mathbb R^{N_y}}\BS G(\BS{\bar x}+\BS q,\BS y)\cdot
\nabla_{\BS y}v_2\dif\mu_{\BS{\bar x}+\BS q}(\BS y)+\\&+
\int_{\mathbb R^{N_y}}\big[\BS F(\BS{\bar x}+\BS q,\BS y)-\BS{\bar F}(\BS{\bar
x})\big]\cdot\nabla_{\BS q}v_1\dif\mu_{\BS{\bar x}+\BS q}(\BS y).
\end{split}
\end{equation}
For the first integral above we express the gradient term as a time
derivative of the function along the flow (in the same way as above
for $v_1$):
\begin{equation}
\BS G(\BS{\bar x}+\BS q,\BS y)\cdot\nabla_{\BS y}v_2 (\tau,\BS{\bar
  x},\BS q,\BS y)=\left.\parderiv{}s v_2(\tau,\BS{\bar x},\BS
q,\phi_{\BS{\bar x}+\BS q} ^s\BS y)\right|_{s=0}.
\end{equation}
The invariant distribution measure $\mu_{\BS{\bar x}+\BS q}$ preserves the
averages of functions of $\phi_{\BS{\bar x}+\BS q}^s\BS y$:
\begin{equation}
\int_{\mathbb R^{N_y}}v_2(\tau,\BS q,\phi_{\BS{\bar x}+\BS q}^s\BS y)
\dif\mu_{\BS{\bar x}+\BS q}(\BS y)=
\mbox{ constant, for all }s.
\end{equation}
This, in turn, results in
\begin{equation}
\parderiv{}s\int_{\mathbb R^{N_y}}v_2(\tau,\BS{\bar x},\BS q,\phi_{\BS{\bar x}+\BS q}^s\BS y)
\dif\mu_{\BS{\bar x}+\BS q}(\BS y)=0,\mbox{ for all $s$, including }s=0,
\end{equation}
and, therefore,
\begin{equation}
\parderiv{v_0}\tau=\int_{\mathbb R^{N_y}}\big[\BS F(\BS{\bar x}+\BS q,\BS y)-\BS{\bar F}
(\BS{\bar x})\big]\cdot\nabla_{\BS q}v_1\dif\mu_{\BS{\bar x}+\BS q}(\BS y).
\end{equation}
\end{itemize}
At this point, substituting the expression for $v_1$ from
\eqref{eq:v11}, and rescaling time $\tau$ back to $t$ yields
\begin{subequations}
\begin{equation}
\parderiv{v_0}t=\left[\BS Q(\BS{\bar x}+\BS q,\BS{\bar x})\cdot\nabla_{\BS q}
+\frac 12\BS S(\BS{\bar x}+\BS q,\BS{\bar x}):
\big(\nabla_{\BS q}\otimes\nabla_{\BS q}\big)\right]v_0,
\end{equation}
\begin{equation}
\BS Q(\BS a,\BS b)=\varepsilon\int_0^\infty\int_{\mathbb R^{N_y}}
\left(\parderiv{}{\BS a}\BS F(\BS a,\phi_{\BS a}^s\BS y)\right)
\big[\BS F(\BS a,\BS y)-\BS{\bar F}(\BS b)\big]
\dif\mu_{\BS a}(\BS y)\dif s,
\end{equation}
\begin{equation}
\label{eq:S}
\BS S(\BS a,\BS b)=2\varepsilon\sym\int_0^\infty\int_{\mathbb R^{N_y}}\big[\BS
  F(\BS a,\phi_{\BS a}^s\BS y)-\BS{\bar F}
(\BS b)\big]\otimes\big[\BS F(\BS a,\BS y)-
\BS{\bar F}(\BS b)\big]\dif\mu_{\BS a}(\BS y)\dif s,
\end{equation}
\end{subequations}
where ``$:$'' denotes the Frobenius product of two matrices, and
``$\sym$'' denotes the symmetric part of a matrix (the skew-symmetric
part is canceled out by the Frobenius product with a symmetric
matrix).

As a result, the next-order correction $\BS{\bar q}$ to the averaged
dynamics in \eqref{eq:dyn_sys_slow_limiting_x}, generated by the
backward Kolmogorov equation above, is given by the It\^o diffusion
process
\begin{subequations}
\label{eq:dx_dq}
\begin{equation}
\deriv{\BS{\bar x}}t=\BS{\bar F}(\BS{\bar x}),
\end{equation}
\begin{equation}
\dif\BS{\bar q}=\BS Q(\BS{\bar x}+\BS{\bar q},\BS{\bar x})\dif t
+\BS\sigma(\BS{\bar x}+\BS{\bar q},\BS{\bar x})\dif\BS W_t,
\end{equation}
\end{subequations}
as long as $\BS{\bar q}$ is small enough. Above, $\BS W_t$ is a
$K$-dimensional Wiener process for some positive integer $K$, $\BS
Q(\BS{\bar x}+\BS{\bar q},\BS{\bar x})$ is the $O(\varepsilon)$ $N_y$-vector drift term, and
$\BS\sigma(\BS{\bar x}+\BS{\bar q},\BS{\bar x})$ is the $O(\sqrt\varepsilon)$ $N_y\times K$
stochastic diffusion matrix, given (non-uniquely) by
\begin{equation}
\label{eq:sigma_S}
\BS\sigma\BS\sigma^T=\BS S.
\end{equation}
Thus far, in \eqref{eq:dx_dq} we obtained the equations for two sets
of variables, $\BS{\bar x}$ and $\BS{\bar q}$, while the evolution of
the stochastic reduced dynamics for \eqref{eq:dyn_sys} is given by the
sum $(\BS{\bar x}+\BS{\bar q})$. However, what we actually need is an
approximate reduced model for $(\BS{\bar x}+\BS{\bar q})$
directly. For that, we add the two equations in \eqref{eq:dx_dq} to
obtain
\begin{equation}
\label{eq:x+q}
\begin{split}
\dif(\BS{\bar x}+\BS{\bar q})&=\BS{\bar F}(\BS{\bar x}+\BS{\bar q})\dif
t +\BS\sigma(\BS{\bar x}+\BS{\bar q},\BS{\bar x}+\BS{\bar q})\dif\BS W_t+\\
+&\big[\BS{\bar F}(\BS{\bar x}) -\BS{\bar F}(\BS{\bar x}+\BS{\bar q})+\BS
  Q(\BS{\bar x}+\BS{\bar q},\BS{\bar x})\big]\dif t+\\
&+[\BS\sigma(\BS{\bar x}+\BS{\bar q},\BS{\bar x})-
\BS\sigma(\BS{\bar x}+\BS{\bar q},\BS{\bar x}+\BS{\bar q})]\dif\BS W_t.
\end{split}
\end{equation}
Above, observe that the first term in the right-hand side is $O(1)$,
the second term is $O(\sqrt\varepsilon)$, and the rest of the terms
are $O(\varepsilon)$ or higher (taking into account that $\BS{\bar q}$
is $O(\varepsilon)$). Additionally, the terms in the second and third
line of \eqref{eq:x+q} are very difficult to compute in practice for
complex nonlinear dynamics. Therefore, we delete the terms in the
second and third lines from \eqref{eq:x+q}, thus completely decoupling
\eqref{eq:x+q} from \eqref{eq:dx_dq}. Then, after replacing $(\BS{\bar
  x}+\BS{\bar q})\to\BS{\bar x}$, we obtain the stochastic reduced
model in \eqref{eq:dyn_sys_slow_limiting_stoch}, where the diffusion
matrix $\BS\sigma$ is computed according to \eqref{eq:S} and
\eqref{eq:sigma_S} with $\BS a=\BS b=\BS{\bar x}$ (further we denote
$\BS S(\BS{\bar x})\equiv\BS S(\BS{\bar x},\BS{\bar x})$).

\subsection*{Remark: $\BS\sigma(\BS{\bar x})$ does not depend on
$\varepsilon$} It is important to note that $\BS S(\BS{\bar x})$, and,
therefore, $\BS\sigma(\BS{\bar x})$, do not in fact depend on
$\varepsilon$. Indeed, observe that
\begin{equation}
\begin{split}
\BS S(\BS{\bar x})=2\varepsilon\sym\int_0^\infty\int_{\mathbb R^{N_y}}\big[\BS F(\BS{\bar
    x},\phi_{\BS{\bar x}}^s\BS y)-\BS{\bar F}(\BS{\bar x})\big] \otimes\big[\BS
  F(\BS{\bar x},\BS y)-\BS{\bar F}(\BS{\bar x})\big] \dif\mu_{\BS{\bar
    x}}(\BS y)\dif s=\\=2\sym\int_0^\infty\int_{\mathbb
  R^{N_y}}\big[\BS F(\BS{\bar x},\phi_{\BS{\bar x}}^{s/\varepsilon}\BS y)-\BS{\bar
    F}(\BS{\bar x})\big] \otimes\big[\BS F(\BS{\bar x},\BS y)-\BS{\bar
    F}(\BS{\bar x})\big] \dif\mu_{\BS{\bar x}}(\BS y)\dif s,
\end{split}
\end{equation}
where $\phi_{\BS{\bar x}}^{s/\varepsilon}\BS y$ is the solution of the
fast dynamics from \eqref{eq:dyn_sys} (with $\varepsilon^{-1}$ still
in front of $\BS G(\BS x,\BS y)$) with $\BS x=\BS{\bar x}$ fixed as a
constant parameter. Thus, there is no need to know what $\varepsilon$
is to compute $\BS S(\BS{\bar x})$, and, in fact, $\varepsilon$ does
not have to be an explicit scaling parameter in the multiscale
dynamics in \eqref{eq:dyn_sys}. This makes the new method practical in
realistic applications, where time-scale difference might not be
explicitly available in the form of a parameter.

\section{Practical implementation of the reduced stochastic model
for a general multiscale process with linear coupling}
\label{sec:practical}

As formulated above in Section \ref{sec:theory}, the new method is
applicable for a broad range of dynamical systems with general forms
of coupling. However, for the clarity of presentation, we dedicate a
separate section to the implementation of the developed method for a
multiscale process with linear coupling. The linear coupling is the
most basic form of coupling in physical processes, however, because of
that it is also probably the most common form of coupling. Below we
describe the complete step-by-step assembly of the reduced model, with
both the deterministic and stochastic terms which parameterize
coupling.

Here we consider the special setting of \eqref{eq:dyn_sys} with linear
coupling between $\BS x$ and $\BS y$:
\begin{equation}
\label{eq:dyn_sys_linear}
\deriv{\BS x}t=\BS f(\BS x)+\BS L_y\BS y,\qquad
\deriv{\BS y}t=\BS g(\BS y)+\BS L_x\BS x,
\end{equation}
where $\BS f$ and $\BS g$ are nonlinear vector functions of $\BS x$
and $\BS y$, respectively, and $\BS L_x$ and $\BS L_y$ are constant
matrices of suitable sizes. As before, $\BS x$ and $\BS y$ denote the
slow and fast variables, respectively. However, one important
distinction between \eqref{eq:dyn_sys_linear} and \eqref{eq:dyn_sys}
is that we no longer make use of the time scale separation parameter
$\varepsilon$; here the assumption is that it is somehow known that
$\BS x$-variables are slow, and $\BS y$-variables are fast, but no
further information about time scale separation is available beyond
that, and, in particular, no explicit time scale separation parameter
is known. As before, we introduce the corresponding fast limiting
system
\begin{equation}
\label{eq:dyn_sys_fast_limiting_linear_z}
\deriv{\BS z}t=\BS g(\BS z)+\BS L_x\BS x
\end{equation}
with $\BS x$ specified as a constant parameter. It is easy to see that
for the multiscale system with linear coupling in
\eqref{eq:dyn_sys_linear}, the matrix $\BS S(\BS{\bar x})$ from
\eqref{eq:S} is computed as the time-lag correlation matrix
\begin{equation}
\label{eq:S_linear}
\BS S(\BS{\bar x})=2\BS L_y\sym\left[\int_0^\infty\int_{\mathbb R^{N_y}}
\big(\phi_{\BS{\bar x}}^s\BS y-\BS{\bar z}(\BS{\bar x})\big)
\big(\BS y-\BS{\bar z}(\BS{\bar x})\big)^T
\dif\mu_{\BS{\bar x}}(\BS y)\dif s\right]\BS L_y^T,
\end{equation}
where $\phi_{\BS{\bar z}}^s$ is the solution of
\eqref{eq:dyn_sys_fast_limiting_linear_z}, with $\BS x$ fixed at
constant parameter $\BS{\bar x}$, forward in time $s$ from an initial
condition $\BS y$, and $\BS{\bar z}(\BS{\bar x})$ is the mean state of
\eqref{eq:dyn_sys_fast_limiting_linear_z}. Before constructing the
reduced model, we need realistic assumptions on what we can compute in
the multiscale system in \eqref{eq:dyn_sys_linear}, which we take
directly from \cite{Abr9,Abr10}:
\begin{itemize}
\item First, we are going to presume that the multiscale dynamical
  system in \eqref{eq:dyn_sys_linear} is not necessarily computable at
  will for arbitrarily long time intervals. The reason is that if it
  is possible, then the need for a reduced model becomes somewhat
  difficult to justify.
\item Even if the full multiscale model is not computable at will, we
  still need some statistical information about it to formulate the
  reduced model. Here, we presume that some typical state $\BS x^*$ of
  the slow variables $\BS x$ is available, such that the dynamics
  evolve in the proximity of $\BS x^*$. For example, a rough estimate
  of the mean state of the slow variables of the full multiscale
  system can be taken as $\BS x^*$, or a nearby state.
\item We presume that the limiting fast dynamics in
  \eqref{eq:dyn_sys_fast_limiting_linear_z} is computable beyond the
  mixing time scale, so that time averages of
  \eqref{eq:dyn_sys_fast_limiting_linear_z} can be computed, at least
  for a single given value $\BS x=\BS x^*$.
\item The assumptions above are the same as in \cite{Abr9,Abr10}, in
  order to ensure compatibility of the proposed stochastic reduced
  model method with what we have already developed in
  \cite{Abr9,Abr10} for deterministic reduced models of multiscale
  dynamics with nonlinear and multiplicative coupling.
\end{itemize}
Under the above assumptions and the ergodicity hypothesis, we compute
the mean state $\langle\BS z\rangle$ and the time covariance matrix
$\BS C(\tau)$ from the long-term time series of the solution $\BS
z(t)$ of \eqref{eq:dyn_sys_fast_limiting_linear_z} with
fixed parameter $\BS x=\BS x^*$:
\begin{subequations}
\begin{equation}
\langle\BS z\rangle=\lim_{T\to\infty}\frac 1T\int_0^T\BS z(t)\dif t,
\end{equation}
\begin{equation}
\BS C(\tau)=\lim_{T\to\infty}\frac 1T\int_0^T\left(\BS z(t+\tau)
-\langle\BS z\rangle\right)\left(\BS z(t)
-\langle\BS z\rangle\right)^T\dif t,
\end{equation}
\begin{equation}
\BS{\bar C}=\int_0^\infty\BS C(\tau)\dif\tau,
\end{equation}
\end{subequations}
Then, we assemble the reduced stochastic model for
\eqref{eq:dyn_sys_linear} in the explicit form as
\begin{equation}
\label{eq:dyn_sys_slow_limiting_linear_stoch}
\dif\BS{\bar x}=\left[\BS f(\BS{\bar x})+\BS L_y\langle\BS z\rangle
+\BS L_y\BS R\BS L_x(\BS{\bar x}-\BS x^*)\right]\dif t+\BS\sigma
\dif\BS W_t,
\end{equation}
where the constant terms $\BS R$ and $\BS\sigma$
are computed as
\begin{subequations}
\begin{equation}
\BS R=\BS{\bar C}\BS C^{-1}(0),
\end{equation}
\begin{equation}
\label{eq:sigma_model}
\BS\sigma\BS\sigma^T=\BS S=
\BS L_y\left(\BS{\bar C}
+\BS{\bar C}^T\right)\BS L_y^T.
\end{equation}
\end{subequations}
Above, the deterministic part in
\eqref{eq:dyn_sys_slow_limiting_linear_stoch} is computed as described
in \cite{Abr9} using the quasi-Gaussian linear response approximation
\cite{MajAbrGro}, while the stochastic part is computed according to
\eqref{eq:S_linear} (where the measure average is replaced with the
time average) with $\BS{\bar x}=\BS x^*$. Observe that the matrix
$\BS\sigma$ will be computed only once, for the particular value $\BS
x^*$. That, in effect, makes it a constant matrix approximation to the
exact $\BS\sigma(\BS{\bar x})$ from \eqref{eq:sigma_S}, provided that
the trajectory $\BS{\bar x}(t)$ of the reduced model in
\eqref{eq:dyn_sys_slow_limiting_linear_stoch} is in the vicinity of
$\BS x^*$. This type of stochastic parameterization with a constant
diffusion matrix is called the additive noise parameterization (as
opposed to the multiplicative noise parameterization, where
$\BS\sigma$ is a function of $\BS{\bar x}$). In the future work, we
plan to extend the method described here onto the multiplicative noise
parameterization.

\subsection*{The choice of $\BS\sigma$}

Observe that $\BS\sigma$ is not determined uniquely by
\eqref{eq:sigma_model}, as multiple decompositions of $\BS S$ into
$\BS\sigma\BS\sigma^T$ are available (the Cholesky decomposition into
the product of a lower-triangular matrix with its own transpose being
one of the examples). Here we use the following algorithm to compute
$\BS\sigma$: first, solve the eigenvalue problem
\begin{equation}
\BS S\BS X=\BS X\BS\Lambda,
\end{equation}
where $\BS X$ is the matrix of eigenvectors, and $\BS\Lambda$ is the
diagonal matrix of eigenvalues. Since $\BS S$ is symmetric and
nonnegative-definite \cite{PavStu}, then all eigenvalues in
$\BS\Lambda$ are guaranteed to be positive, and all eigenvectors in
$\BS X$ are guaranteed to be orthonormal (so that $\BS X^{-1}=\BS
X^T$). Then, we compute $\BS\sigma$ as
\begin{equation}
\BS\sigma=\BS X\BS\Lambda^{1/2}\BS X^T.
\end{equation}
This method uniquely determines $\BS\sigma$ as a square, symmetric,
and positive-definite diffusion matrix for the Wiener process $\BS
W_t$.

\section{Computational study: the two-scale Lorenz 96 model with linear coupling}
\label{sec:numerical}

The test multiscale system used to study the new method of stochastic
parameterization here is the rescaled two-scale Lorenz 96 system with
linear coupling, previously used in \cite{Abr9} to study the
deterministic reduced model parameterization. The Lorenz 96 system is
given by
\begin{subequations}
\label{eq:lorenz_two_scale}
\begin{equation}
\dot x_i=x_{i-1}(x_{i+1}-x_{i-2})+\frac 1{\beta_x}(\bar x(x_{i+1}-x_{i-2})
-x_i)+\frac{F_x-\bar x}{\beta_x^2}-\frac{\lambda_y}J\sum_{j=1}^Jy_{i,j},
\end{equation}
\begin{equation}
\dot y_{i,j}=\frac 1\varepsilon\bigg[y_{i,j+1}
(y_{i,j-1}-y_{i,j+2})+\frac 1{\beta_y}(\bar y(y_{i,j-1}-y_{i,j+2})-y_{i,j})
+\frac{F_y-\bar y}{\beta_y^2}\bigg]+\frac{\lambda_x}\varepsilon x_i,
\end{equation}
\end{subequations}
with $1\leq i\leq N_x$ and $1\leq j\leq J$, such that $N_y=N_xJ$. The
model has periodic boundary conditions: $x_{i+N_x}=x_i$,
$y_{i,j+J}=y_{i+1,j}$, and $y_{i+N_x,j}=y_{i,j}$. The parameter
$\varepsilon\ll 1$ sets the time scale separation between the slow
variables $x_i$, and the fast variables $y_{i,j}$. The parameters
$(\bar x,\beta_x)$ and $(\bar y,\beta_y)$ are the (mean, standard
deviation) pairs for the corresponding uncoupled and unrescaled Lorenz
models
\begin{subequations}
\begin{equation}
\dot x_i=x_{i-1}(x_{i+1}-x_{i-2})-x_i+F_x,
\end{equation}
\begin{equation}
\dot y_{i,j}=y_{i,j+1}(y_{i,j-1}-y_{i,j+2})-y_{i,j}+F_y,
\end{equation}
\end{subequations}
with the same periodic boundary conditions. The rescaling above
ensures that the Lorenz 96 model in \eqref{eq:lorenz_two_scale} has
zero mean state and unit standard deviation for both slow and fast
variables in the absence of coupling ($\lambda_x=\lambda_y=0$), and
remain near these values when $\lambda_x$ and $\lambda_y$ are nonzero.
The linear coupling above preserves the energy of the form
\begin{equation}
E=\frac{\lambda_x}2\sum_{i=1}^{N_x}x_i^2+\frac{\varepsilon\lambda_y}{2J}
\sum_{i=1}^{N_x}\sum_{j=1}^Jy_{i,j}^2.
\end{equation}

Below we display the results of a numerical study of the stochastic
reduced model for slow variables of multiscale dynamics, using the
rescaled Lorenz system in \eqref{eq:lorenz_two_scale} as the test
model. We compare the statistical properties of the slow variables for
the four following systems:
\begin{enumerate}
\item The complete rescaled Lorenz system from
  \eqref{eq:lorenz_two_scale};
\item The stochastic reduced model from
  \eqref{eq:dyn_sys_slow_limiting_linear_stoch};
\item The deterministic reduced model obtained from
  \eqref{eq:dyn_sys_slow_limiting_linear_stoch} by removing the
  stochastic forcing (this model was previously developed in
  \cite{Abr9});
\item The poor man's version of
  \eqref{eq:dyn_sys_slow_limiting_linear_stoch} with no stochastic
  forcing and the first-order linear deterministic correction term
  $\BS R$ set to zero (further referred to as the ``zero-order''
  system, same as in \cite{Abr9}). This zero-order system represents
  the simplest reduced model with constant parameterization of
  coupling terms.
\end{enumerate}
The fixed parameter $\BS x^*$ for the computation of $\langle\BS
z\rangle$, $\BS R$ and $\BS\sigma$ was set to the long-term mean state
$\langle\BS x\rangle$ of the full multiscale model in
\eqref{eq:lorenz_two_scale} (in computationally expensive models, a
rough estimate could be used). Like in \cite{Abr9,Abr10}, for all
computational statistical results presented below, the averaging time
window $T_{av}$ equals 10000 time units.

Due to translational invariance of the studied models, the statistics
are invariant with respect to the index shift for the variables
$x_i$. For diagnostics, we monitor the following long-term statistical
quantities of $x_i$:
\begin{enumerate}[a.]
\item The distribution density functions, computed by bin-counting. A
  distribution density function gives the best information about the
  one-point statistics of $x_i$, as it shows the statistical
  distribution of $x_i$ in the phase space.
\item The time auto-correlation functions $\langle
  x_i(t)x_i(t+s)\rangle$, where the time average is over $t$,
  normalized by the variance $\langle x_i^2\rangle$ (so that it always
  starts with 1).
\item The time cross-correlation functions $\langle
  x_i(t)x_{i+1}(t+s)\rangle$, also normalized by the variance $\langle
  x_i^2\rangle$.
\item The energy auto-correlation function $$K(s)=\frac{\langle
  x_i^2(t)x_i^2(t+s)\rangle}{\langle x_i^2\rangle^2+2\langle
  x_i(t)x_i(t+s)\rangle^2}.$$ This energy auto-correlation function
  measures the non-Gaussianity of the process (it is identically 1 for
  all $s$ if the process is Gaussian, such as the Ornstein-Uhlenbeck
  process). For details, see \cite{MajTimVan4}.
\end{enumerate}
The following dynamical regimes are studied:
\begin{itemize}
\item $N_x=20$, $J=4$ (so that $N_y=80$). Thus, the number of the fast
  variables is four times greater than the number of the slow variables.
\item $\varepsilon=0.01,0.1$. The time scale separation of two orders
  of magnitude ($\varepsilon=0.01$) is consistent with typical
  real-world geophysical processes (for example, the annual and
  diurnal cycles of the Earth's atmosphere). Also, in large-scale
  atmospheric processes the time scale separation between slow and
  fast variables can be weaker than that (for example, typical time
  scale of equatorial Kelvin waves is about 70 days, and that of Yanai
  waves is about 30 days, which is well in between the annual and
  diurnal cycles), so we additionally test the dynamical regimes with
  weak time scale separation $\varepsilon=0.1$.
\item $\lambda_x=\lambda_y=0.3,0.35$. These values of coupling are
  chosen so that they are neither too weak, nor too strong (although
  0.3 is weaker, and 0.35 is stronger). However, this small variation
  in coupling changes the dynamical regime in the slow variables
  between moderately ($\varepsilon=0.3$) and weakly
  ($\varepsilon=0.35$) chaotic and mixing, due to the suppression of
  chaos effect previously studied in \cite{Abr8}.
\item $F_x=6$. The slow forcing $F_x$ adjusts the chaos and mixing
  properties of the slow variables, and in this work it is set to a
  weak-to-moderate chaotic regime $F_x=6$. The reason for that is that
  it was found previously in \cite{Abr9} that in strongly chaotic and
  mixing regimes at slow variables there is not much of a difference
  between the multiscale dynamics and reduced models.
\item $F_y=16$. The fast forcing adjusts the chaos and mixing
  properties of the fast variables. Here the value of $F_y$ is chosen
  so that the fast variables are strongly chaotic and mixing for
  $F_y=16$.
\end{itemize}

\subsection{Moderate mixing at slow variables with weak time scale separation}

\begin{figure}
\picturehere{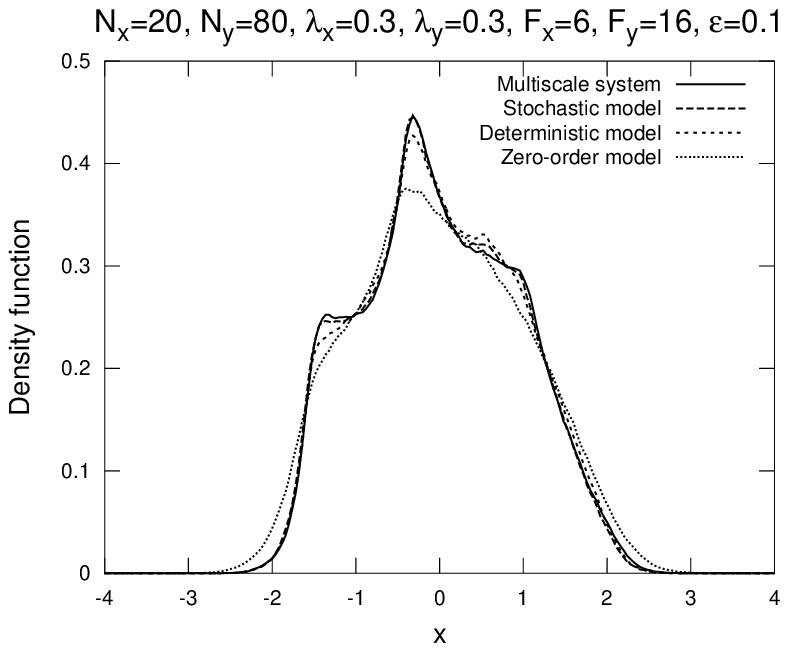}%
\picturehere{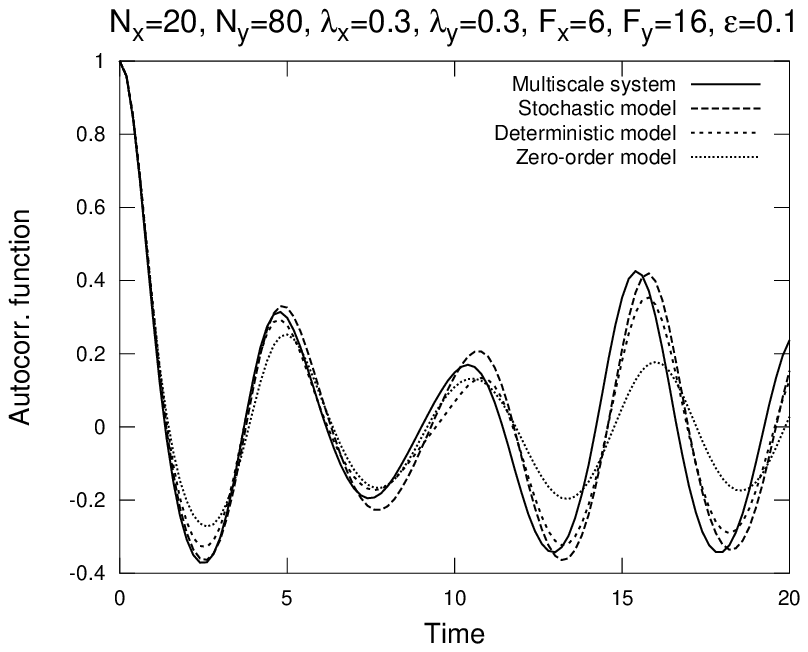}\\%
\picturehere{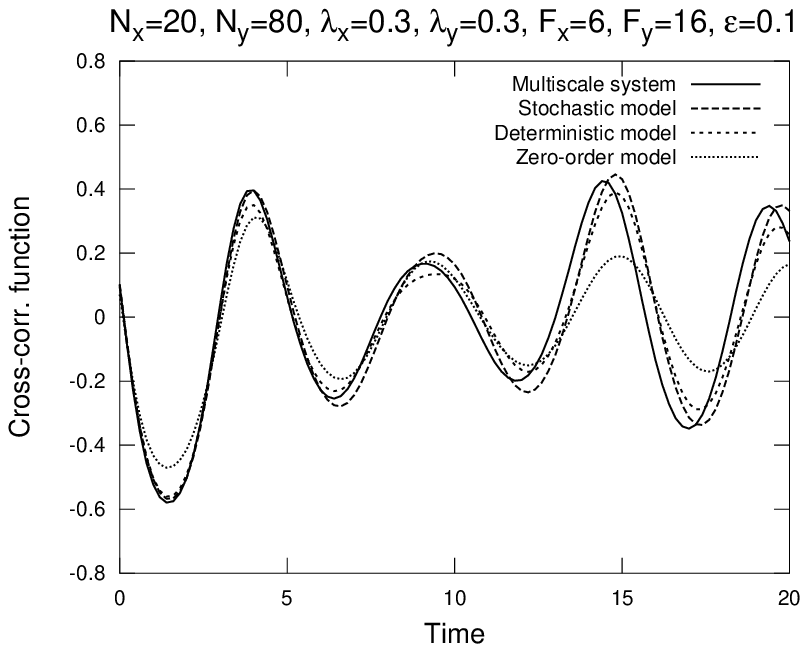}%
\picturehere{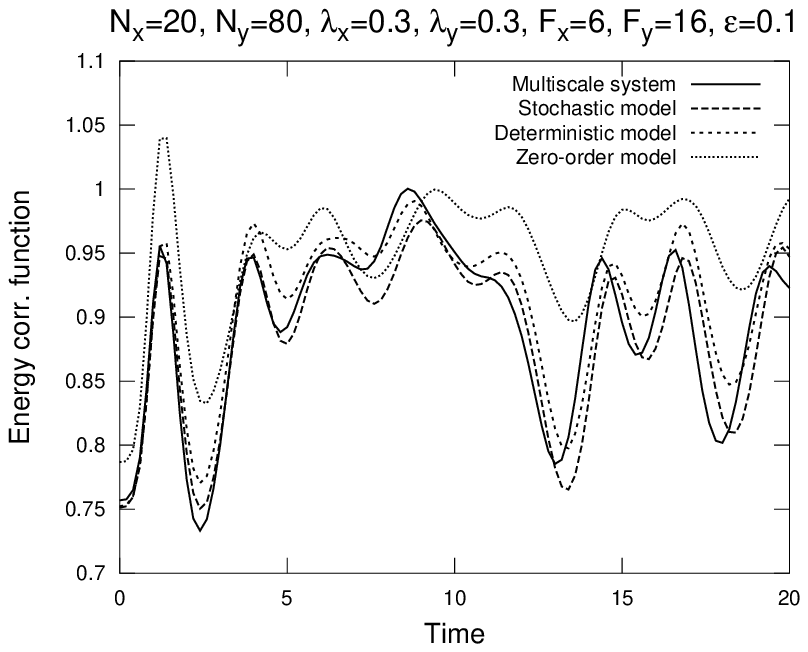}%
\caption{Upper-left -- distribution density function, upper-right --
  time auto-correlation correlation function, lower-left -- time
  cross-correlation function, lower-right -- energy auto-correlation
  function. $F_x=6$, $F_y=16$, $\lambda_x=\lambda_y=0.3$,
  $\varepsilon=0.1$.}
\label{fig:Fx6_Fy16_lx0.3_ly0.3_eps0.1}
\end{figure}

Here we present the comparison of different statistics for the
dynamical regime with moderate chaos and mixing at slow variables
(achieved by setting $\lambda_x=\lambda_y=0.3$) and weak time scale
separation (achieved by setting $\varepsilon=0.1$). In Figure
\ref{fig:Fx6_Fy16_lx0.3_ly0.3_eps0.1} we show the distribution density
functions, time auto-correlation functions, time cross-correlation
functions, and energy auto-correlation functions for the multiscale
dynamics in \eqref{eq:lorenz_two_scale}, and three different kinds of
the reduced models: the new stochastic reduced model, the
deterministic reduced model from \cite{Abr9}, and the zero-order
reduced model with constant parameterization of coupling terms (which
is given for reference as a simplest/poorest form of
parameterization). Observe that the difference between the new
stochastic and deterministic models here are rather small, however, it
does look like the additional stochastic term improves the
statistics. In particular, the stochastic terms makes the distribution
density more spiky towards the multiscale density, and appears to
produce better fits for correlation functions. The probable reason why
there is not much of a difference between the deterministic and
stochastic reduced models is that the slow dynamics are not
sufficiently weakly chaotic for the stochastic term to make a
significant difference. Table \ref{tab:Fx6_Fy16_lx0.3_ly0.3_eps0.1}
confirms that there is some improvement error-wise in the distribution
density and energy auto-correlation, but no improvement in auto- and
cross-correlations (probably due to oscillations getting out of sync
with increased lag time).
\begin{table}
\begin{center}%
\begin{tabular}{|c||c|c|c|}%
\hline%
& Stochastic & Deterministic & Zero-order \\%
\hline\hline%
Density & $3.803\cdot 10^{-3}$ & $7.424\cdot 10^{-3}$ & $2.093\cdot 10^{-2}$ \\%
Corr. & $0.1218$ & $0.1152$ & $0.1935$ \\%
Cross-corr. & $0.1297$ & $0.1222$ & $0.2118$ \\%
Energy corr. & $1.312\cdot 10^{-2}$ & $1.436\cdot 10^{-2}$ & $3.473\cdot 10^{-2}$ \\%
\hline%
\end{tabular}%
\end{center}%
\caption{Relative errors between the slow variables of the full
  multiscale system, and different reduced models, computed for the
  plots in Figure \ref{fig:Fx6_Fy16_lx0.3_ly0.3_eps0.1}. $F_x=6$,
  $F_y=16$, $\lambda_x=\lambda_y=0.3$, $\varepsilon=0.1$.}%
\label{tab:Fx6_Fy16_lx0.3_ly0.3_eps0.1}
\end{table}

\subsection{Moderate mixing at slow variables with strong time scale separation}

\begin{figure}
\picturehere{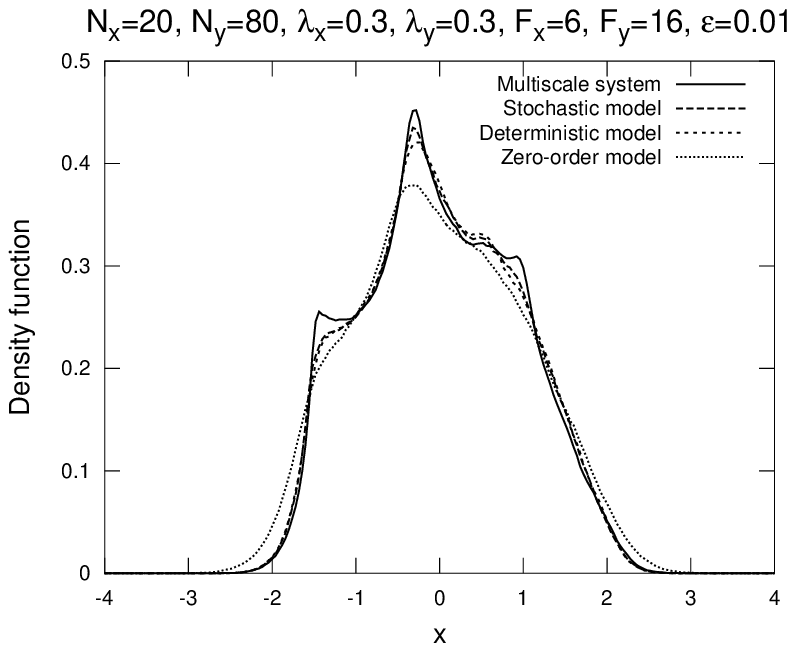}%
\picturehere{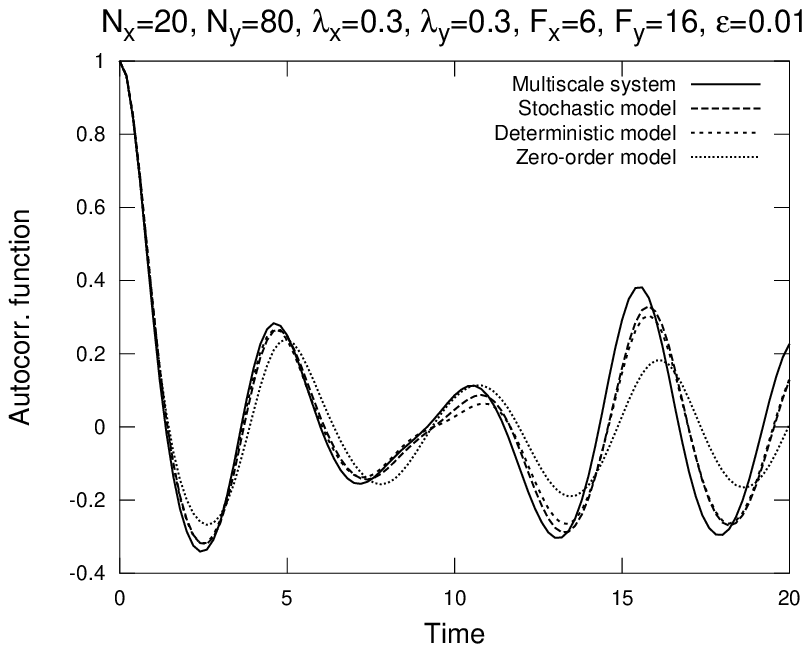}\\%
\picturehere{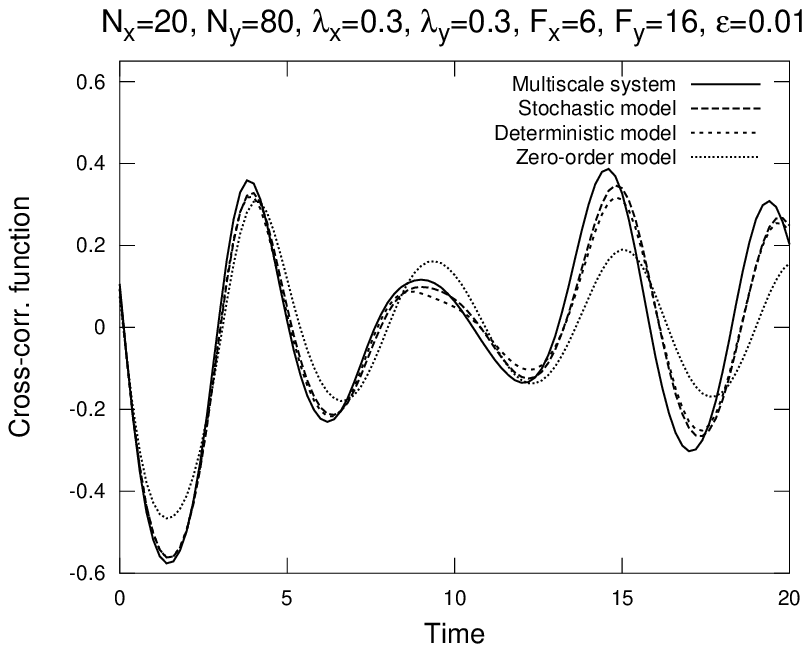}%
\picturehere{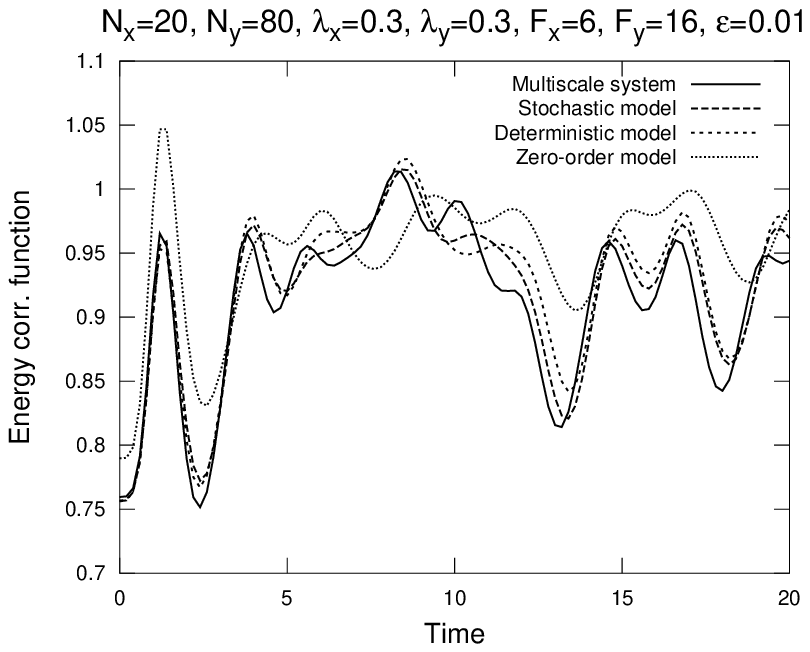}%
\caption{Upper-left -- distribution density function, upper-right --
  time auto-correlation correlation function, lower-left -- time
  cross-correlation function, lower-right -- energy auto-correlation
  function. $F_x=6$, $F_y=16$, $\lambda_x=\lambda_y=0.3$,
  $\varepsilon=0.01$.}
\label{fig:Fx6_Fy16_lx0.3_ly0.3_eps0.01}
\end{figure}

Here we present the comparison of different statistics for the
dynamical regime with moderate chaos and mixing at slow variables
(achieved by setting $\lambda_x=\lambda_y=0.3$) and strong time scale
separation (achieved by setting $\varepsilon=0.01$). In Figure
\ref{fig:Fx6_Fy16_lx0.3_ly0.3_eps0.01} we show the distribution
density functions, time auto-correlation functions, time
cross-correlation functions, and energy auto-correlation functions for
the multiscale dynamics in \eqref{eq:lorenz_two_scale}, and three
different kinds of the reduced models: the new stochastic reduced
model, the deterministic reduced model from \cite{Abr9}, and the
zero-order reduced model with constant parameterization of coupling
terms. Observe that the difference between the new stochastic and
deterministic models here are even smaller than for the regime with
weak time-scale separation, probably due to the fact that the
contribution of the stochastic term scales as $\sqrt\varepsilon$.
However, again, it looks like the additional stochastic term improves
the statistics somewhat. Table \ref{tab:Fx6_Fy16_lx0.3_ly0.3_eps0.01}
show that there is some improvement error-wise in all presented
statistics, but not enough to be discernible visually in Figure
\ref{fig:Fx6_Fy16_lx0.3_ly0.3_eps0.01}.
\begin{table}
\begin{center}%
\begin{tabular}{|c||c|c|c|}%
\hline%
& Stochastic & Deterministic & Zero-order \\%
\hline\hline%
Density & $8.105\cdot 10^{-3}$ & $1.048\cdot 10^{-2}$ & $2.233\cdot 10^{-2}$ \\%
Corr. & $9.309\cdot 10^{-2}$ & $9.627\cdot 10^{-2}$ & $0.1923$ \\%
Cross-corr. & $9.57\cdot 10^{-2}$ & $9.99\cdot 10^{-2}$ & $0.2129$ \\%
Energy corr. & $9.042\cdot 10^{-3}$ & $1.209\cdot 10^{-2}$ & $2.776\cdot 10^{-2}$ \\%
\hline%
\end{tabular}%
\end{center}%
\caption{Relative errors between the slow variables of the full
  multiscale system, and different reduced models, computed for the
  plots in Figure \ref{fig:Fx6_Fy16_lx0.3_ly0.3_eps0.01}. $F_x=6$,
  $F_y=16$, $\lambda_x=\lambda_y=0.3$, $\varepsilon=0.01$.}%
\label{tab:Fx6_Fy16_lx0.3_ly0.3_eps0.01}
\end{table}

\subsection{Weak mixing at slow variables with weak time scale separation}

\begin{figure}
\picturehere{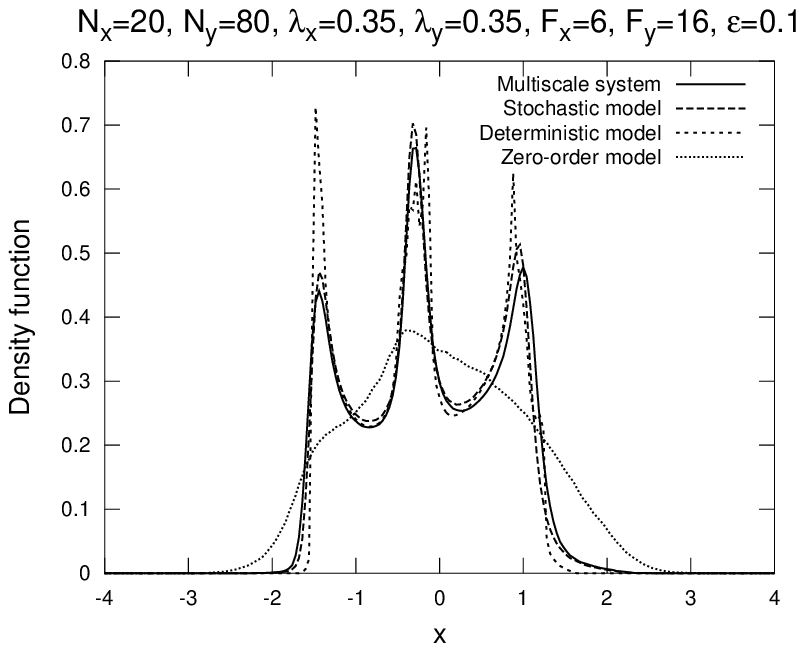}%
\picturehere{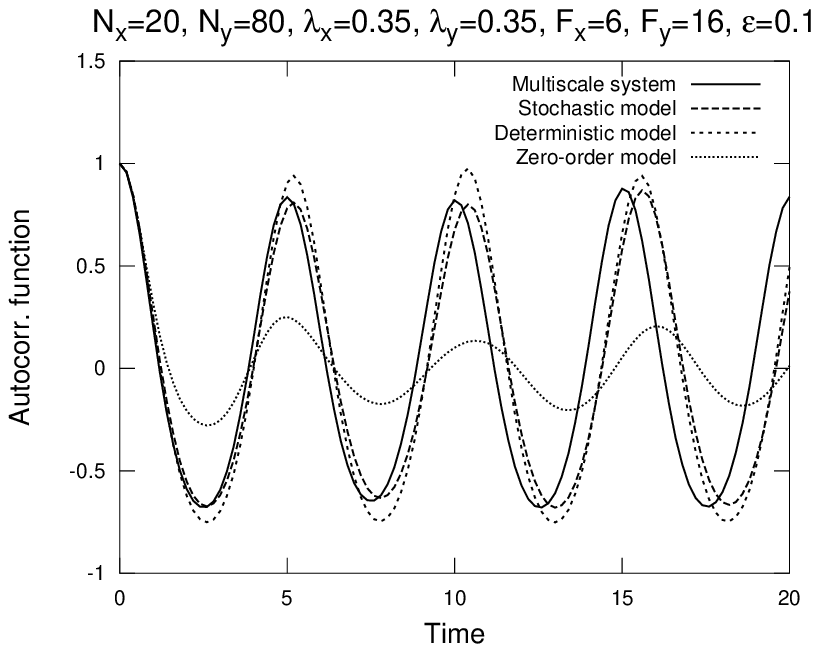}\\%
\picturehere{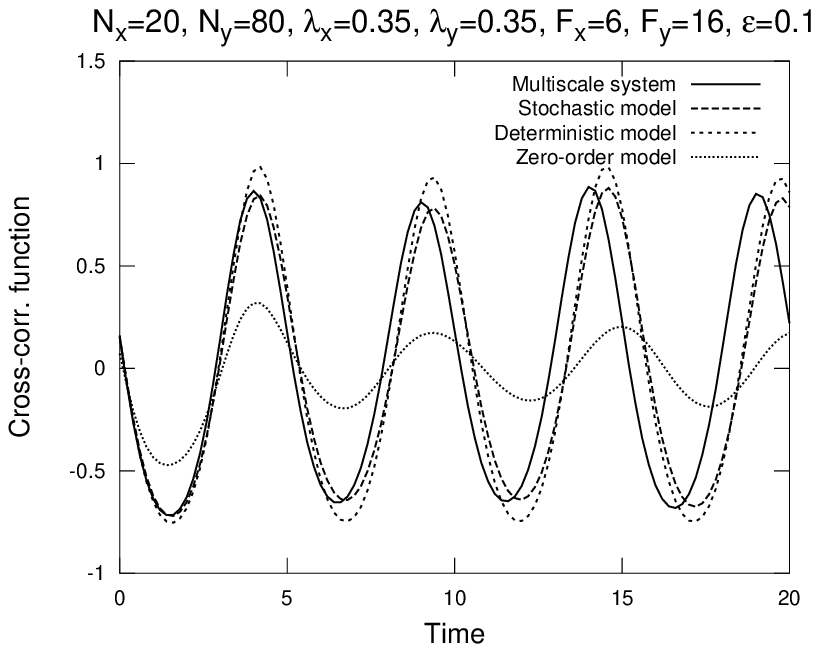}%
\picturehere{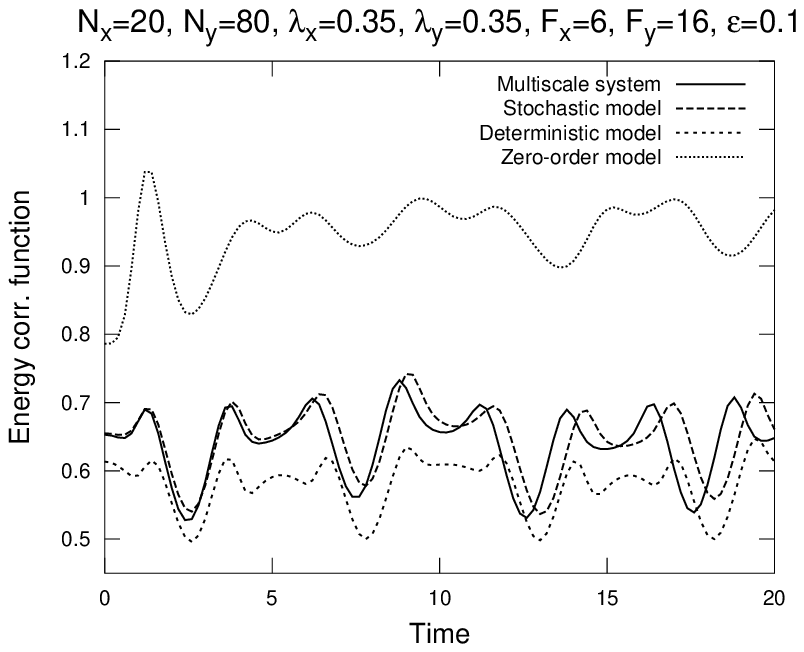}%
\caption{Upper-left -- distribution density function, upper-right --
  time auto-correlation correlation function, lower-left -- time
  cross-correlation function, lower-right -- energy auto-correlation
  function. $F_x=6$, $F_y=16$, $\lambda_x=\lambda_y=0.35$,
  $\varepsilon=0.1$.}
\label{fig:Fx6_Fy16_lx0.35_ly0.35_eps0.1}
\end{figure}
\begin{table}
\begin{center}%
\begin{tabular}{|c||c|c|c|}%
\hline%
& Stochastic & Deterministic & Zero-order \\%
\hline\hline%
Density & $2.166\cdot 10^{-2}$ & $4.83\cdot 10^{-2}$ & $7.516\cdot 10^{-2}$ \\%
Corr. & $0.2322$ & $0.2335$ & $0.3584$ \\%
Cross-corr. & $0.2277$ & $0.2346$ & $0.3557$ \\%
Energy corr. & $2.858\cdot 10^{-2}$ & $5.031\cdot 10^{-2}$ & $0.2163$ \\%
\hline%
\end{tabular}%
\end{center}%
\caption{Relative errors between the slow variables of the full
  multiscale system, and different reduced models, computed for the
  plots in Figure \ref{fig:Fx6_Fy16_lx0.35_ly0.35_eps0.1}. $F_x=6$,
  $F_y=16$, $\lambda_x=\lambda_y=0.35$, $\varepsilon=0.1$.}%
\label{tab:Fx6_Fy16_lx0.35_ly0.35_eps0.1}
\end{table}
Here we present the comparison of different statistics for the
dynamical regime with weak chaos and mixing at slow variables
(achieved by setting $\lambda_x=\lambda_y=0.35$) and weak time scale
separation (achieved by setting $\varepsilon=0.1$). In Figure
\ref{fig:Fx6_Fy16_lx0.35_ly0.35_eps0.1} we show the distribution
density functions, time auto-correlation functions, time
cross-correlation functions, and energy auto-correlation functions for
the multiscale dynamics in \eqref{eq:lorenz_two_scale}, and three
different kinds of the reduced models: the new stochastic reduced
model, the deterministic reduced model from \cite{Abr9}, and the
zero-order reduced model with constant parameterization of coupling
terms. Here we can see a significant improvement between the
deterministic reduced model from \cite{Abr9} and the new stochastic
reduced model. In particular, what apparently happens here is that the
deterministic reduced model from \cite{Abr9} turns out to be less
chaotic and mixing than the original multiscale dynamics at slow
variables (observe that the distribution density is very spiky for the
deterministic reduced model, while the time auto- and
cross-correlation functions are less mixing, and the energy
auto-correlation function is more sub-Gaussian than those for the
multiscale dynamics). Then, the introduction of the stochastic term
results in improvement of mixing and sub-Gaussianity, and also
smoothens out the spikes on the distribution density, resulting in
better approximation of the multiscale dynamics. Table
\ref{tab:Fx6_Fy16_lx0.35_ly0.35_eps0.1} confirms that there is
improvement of statistics error-wise with the introduction of the
stochastic term in the new reduced model.

\subsection{Weak mixing at slow variables with strong time scale separation}

\begin{figure}
\picturehere{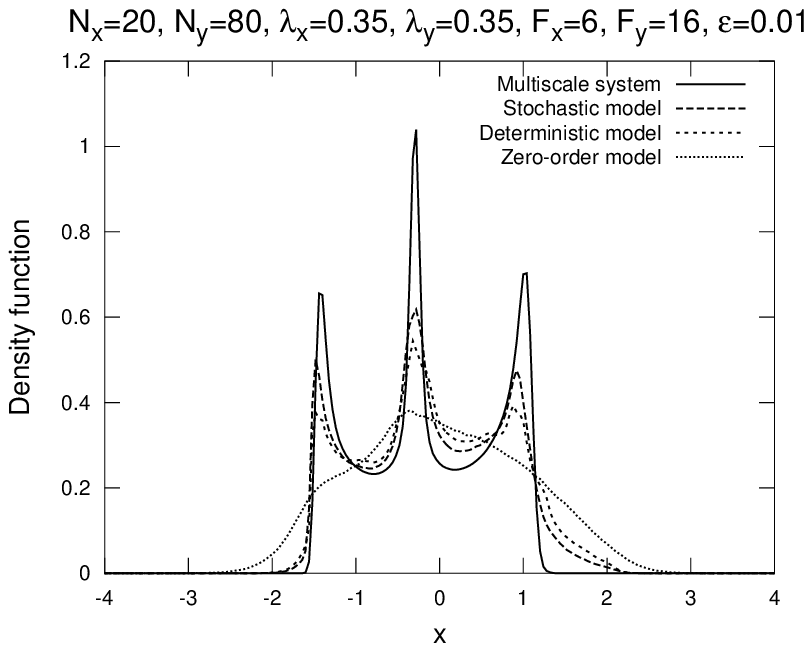}%
\picturehere{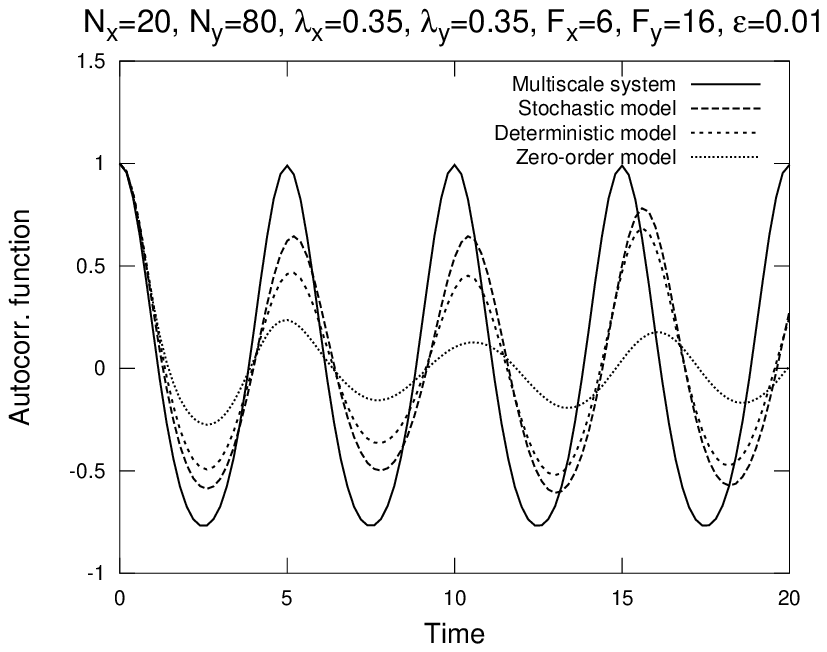}\\%
\picturehere{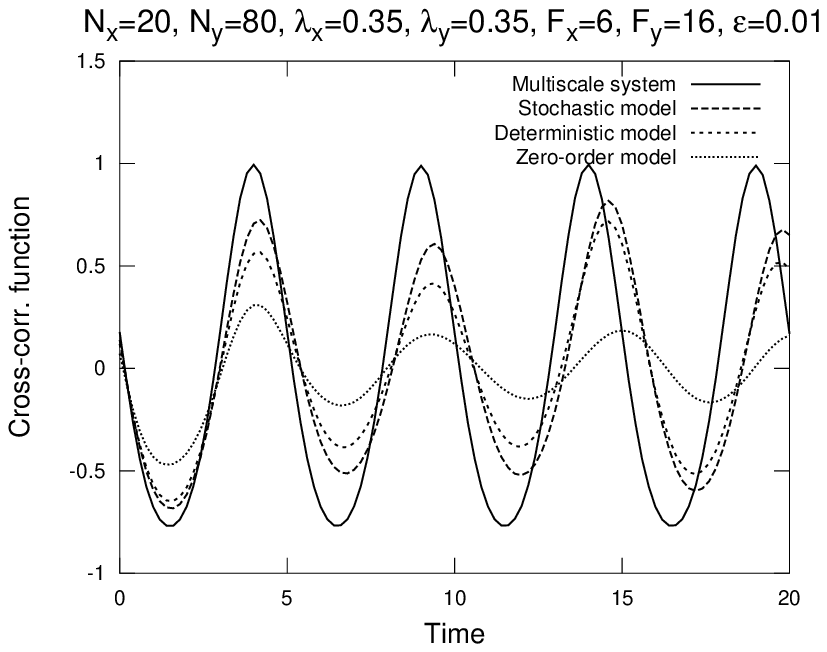}%
\picturehere{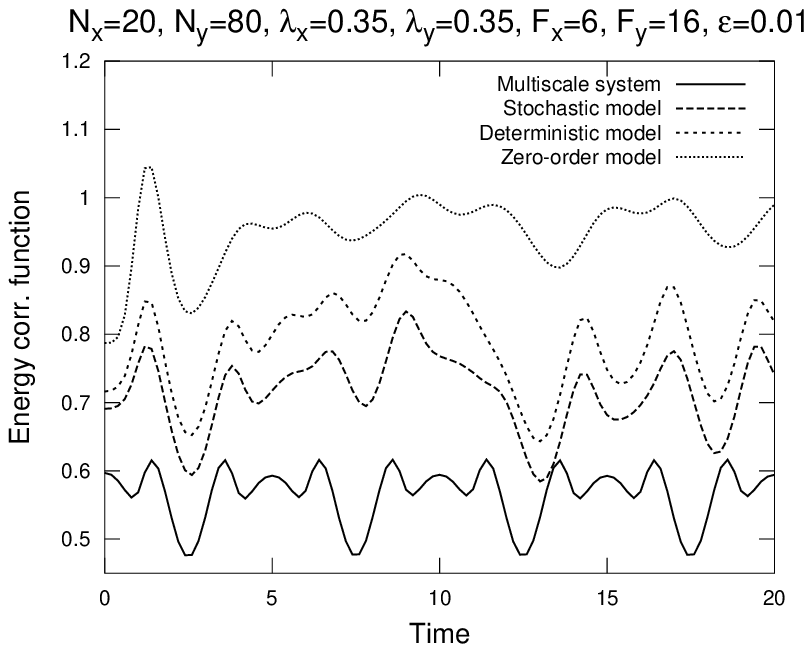}%
\caption{Upper-left -- distribution density function, upper-right --
  time auto-correlation correlation function, lower-left -- time
  cross-correlation function, lower-right -- energy auto-correlation
  function. $F_x=6$, $F_y=16$, $\lambda_x=\lambda_y=0.35$,
  $\varepsilon=0.01$.}
\label{fig:Fx6_Fy16_lx0.35_ly0.35_eps0.01}
\end{figure}
\begin{table}
\begin{center}%
\begin{tabular}{|c||c|c|c|}%
\hline%
& Stochastic & Deterministic & Zero-order \\%
\hline\hline%
Density & $6.237\cdot 10^{-2}$ & $7.716\cdot 10^{-2}$ & $0.1088$ \\%
Corr. & $0.2629$ & $0.2752$ & $0.3769$ \\%
Cross-corr. & $0.2556$ & $0.2684$ & $0.3726$ \\%
Energy corr. & $0.1254$ & $0.1846$ & $0.3059$ \\%
\hline%
\end{tabular}%
\end{center}%
\caption{Relative errors between the slow variables of the full
  multiscale system, and different reduced models, computed for the
  plots in Figure \ref{fig:Fx6_Fy16_lx0.35_ly0.35_eps0.01}. $F_x=6$,
  $F_y=16$, $\lambda_x=\lambda_y=0.35$, $\varepsilon=0.01$.}%
\label{tab:Fx6_Fy16_lx0.35_ly0.35_eps0.01}
\end{table}
Here we present the comparison of different statistics for the
dynamical regime with weak chaos and mixing at slow variables
(achieved by setting $\lambda_x=\lambda_y=0.35$) and strong time scale
separation (achieved by setting $\varepsilon=0.01$). In Figure
\ref{fig:Fx6_Fy16_lx0.35_ly0.35_eps0.01} we show the distribution
density functions, time auto-correlation functions, time
cross-correlation functions, and energy auto-correlation functions for
the multiscale dynamics in \eqref{eq:lorenz_two_scale}, and three
different kinds of the reduced models: the new stochastic reduced
model, the deterministic reduced model from \cite{Abr9}, and the
zero-order reduced model with constant parameterization of coupling
terms. Here, again, we can see a significant improvement between the
deterministic reduced model from \cite{Abr9} and the new stochastic
reduced model, however, in a somewhat reverse way if compared to the
previous set-up with weak time scale separation. Here observe that the
deterministic reduced model is more chaotic and mixing than the slow
variables of the multiscale dynamics, and the stochastic term makes
the new model less chaotic and mixing to better match the multiscale
dynamics. While at first seeming counter-intuitive, this effect can
happen due to the stochastic noise pushing the solution off the
unstable manifold of the deterministic system (where the chaotic and
mixing motion occurs) into the dissipative absorbing region
surrounding the system's attractor while not having enough strength to
compensate for the lack of chaos and mixing with its own random
forcing. In fact, the same (but somewhat weaker) effect can also be
observed in Figures \ref{fig:Fx6_Fy16_lx0.3_ly0.3_eps0.1} and
\ref{fig:Fx6_Fy16_lx0.3_ly0.3_eps0.01} for weaker coupling and
stronger chaos and mixing, so one can presume that it should be
generally common in stochastic reduced models, especially if the time
scale separation is strong and, because of that, the stochastic noise
is weak enough to produce its own mixing. Table
\ref{tab:Fx6_Fy16_lx0.35_ly0.35_eps0.01} confirms that there is some
improvement of statistics error-wise with the introduction of the
stochastic term in the new reduced model.

\section{Conclusions}
\label{sec:conclusions}

In the current work we improve the recently developed method
\cite{Abr9,Abr10} for deterministic reduced models of multiscale
dynamics with the higher-order additive stochastic term which
parameterizes the slow-fast interactions with random noise. The method
is based on the homogenization techniques for multiscale systems
\cite{PavStu} and offers a practical way of creating stochastic
reduced models for a broad range of general multiscale processes. As
demonstrated above in Section \ref{sec:practical}, the stochastic term
upgrade comes at no additional computational cost for systems with
linear coupling between slow and fast variables, as it uses the same
time correlation matrix of the fast variables for a fixed state of the
slow variables as the response term in the deterministic part of
coupling parameterization. Another practical advantage of the new
method is that it does not require an explicit time-scale separation
parameter between the slow and fast variables of multiscale
dynamics. We tested the new method numerically using the two-scale
Lorenz 96 model with linear coupling between the slow and fast
variables in a range of dynamical regimes with weak/strong time-scale
separation, chaos and mixing at the slow variables. The new stochastic
reduced model consistently improved the results of the previously
developed deterministic approach \cite{Abr9,Abr10} in the two
different dynamical regimes:
\begin{enumerate}
\item In the situation where the deterministic reduced model was less
  chaotic and weaker mixing than the slow variables of the full
  multiscale dynamics, the stochastic model was more chaotic and
  stronger mixing, due to stochastic forcing smoothening out spikes in
  distribution density and introducing random decorrelation in time
  auto- and cross-correlation functions.
\item In the situation where the deterministic reduced model was more
  chaotic and stronger mixing than the slow variables of the full
  multiscale dynamics, the stochastic model suppressed chaos and
  mixing in the reduced dynamics. While this effect seems somewhat
  counter-intuitive, it can be explained by the stochastic noise
  pushing the solution off the unstable manifold of the deterministic
  system (where the chaotic and mixing motion occurs) into the
  dissipative region around the system's attractor, while not having
  enough random force to increase chaos and mixing on its own.
\end{enumerate}
In the future work, we plan to extend the new method of stochastic
parameterization of reduced slow dynamics onto multiplicative
stochastic forcing parameterization. Observe that the additive
stochastic forcing parameterization in the current work emerges from
the constant diffusion matrix approximation above in Section
\ref{sec:practical}. However, if the constant diffusion matrix
approximation is improved by including its higher order Taylor
expansion terms (as was done to the deterministic coupling
parameterization in \cite{Abr9,Abr10}), the inclusion of the
higher-order terms will result in the multiplicative coupling with the
diffusion matrix of the reduced model being a function of the slow
variables. For the multiplicative diffusion matrix approximation, we
plan to use an approach similar to that in \cite{Abr9,Abr10}.

{\bf Acknowledgments.} This work was supported by the National Science
Foundation CAREER grant DMS-0845760, and the Office of Naval Research
grants N00014-09-0083 and 25-74200-F6607.

\end{document}